\documentclass [a4paper,10pt] {article}

\usepackage{graphicx}

\usepackage{natbib}

\input xy
\xyoption{all}

\usepackage{amsmath}
\usepackage{amsfonts}
\usepackage{amssymb}

\usepackage[T1]{fontenc}

\DeclareMathOperator\Stab{Stab}
\DeclareMathOperator\tr{tr}

\newtheorem{Def}{Definition}[section]
\newtheorem{Defprop}[Def]{Definition-Proposition}
\newtheorem{Th}{Theorem}
\newtheorem{Prop}{Proposition}[section]
\newtheorem{Lem}[Prop]{Lemma}

\newcommand{\Q}{\mathbb{Q}}

\title{Proof of Stanley's conjecture about irreducible character values of the symmetric group}
\author{Valentin F\'eray}
\date{December 4, 2006}

\begin{document}

\maketitle

\begin{abstract}
In his paper \cite{St1}, R. Stanley finds a nice combinatorial formula for characters of irreducible representations of the symmetric group of rectangular shape. Then, in \cite{St2}, he gives a conjectural generalisation for any shape. Here, we will prove this formula using shifted Schur functions and Jucys-Murphy elements.
\end{abstract}

\section{The main theorem}
In this article, we will think of $S(n)$ as the group of permutations of $\{1,\ldots,n\}$ so there are canonical embeddings of $S_k$ in $S_n$ ($n \geq k$). The decomposition in cycles with disjoint supports of an element $\sigma \in S(n)$ will play a central role, so we will denote by $C(\sigma)$ the corresponding partition of $\{1,\ldots,n\}$.\\
A partition of $n$ is a weakly decreasing sequence of integers of sum $n$. The irreducible representations of $S(n)$ are canonically indexed by partitions $\lambda$ of $n$ (denoted $\lambda \vdash n$) and $\chi_\lambda$ is the notation for the associated character. For $\mu \in S(k)$ and $\lambda \vdash n$, we will look here at the normalised character, defined by :
$$\widehat{\chi}_\lambda(\mu) = \frac{n(n-1)\ldots (n-k+1) \chi_\lambda(\mu)}{\chi_\lambda(Id_n)},$$
where we have to identify $\mu$ with its image by the natural embedding of $S(k)$ in $S(n)$ to compute $\chi_\lambda(\mu)$.\\

Let $m$ be a positive integer.
\begin{itemize} \item If $\textbf{p}=(p_1,\ldots,p_m)$ and $\textbf{q}=(q_1,\ldots,q_m)$ with $q_1 \geq \ldots \geq q_m$ are two sequences of positive integers, we will denote by $\textbf{p} \times \textbf{q}$ the partition :
$$(\textbf{p} \times \textbf{q})_i = \left\{ \begin{array}{lcl}
	q_1 & \text{if} & 1 \leq i \leq p_1; \\ q_2 & \text{if} & p_1+1 \leq i \leq p_1 + p_2;\\ & \vdots & \\ q_m & \text{if} & p_1 + \ldots + p_{m-1}+1 \leq i \leq p_1 + \ldots + p_m;\\0 & \text{if} & i > p_1 + \ldots + p_m.
                                            \end{array} \right.$$
\item Take an element of $S(k)$ and assign to each of its cycle an integer between $1$ and $m$. The set of such coloured permutations (formally of pairs $(\sigma,\varphi)$, $\sigma \in S(k)$ and $\varphi : C(\sigma) \rightarrow \{1,\ldots,m\}$) is denoted by $S(k)^{(m)}$. Note that $\varphi$ can also be seen as a function $\{1,\ldots,k\} \rightarrow \{1,\ldots,m\}$ invariant by action of $\sigma$. Given a coloured permutation $(\sigma,\varphi) \in S(k)^{(m)}$ and a non-coloured one $\mu \in S(k)$, we can define their product (but it doesn't define a right group action) by $ (\sigma,\varphi) \cdot \mu = (\sigma \mu,\psi)$ where,
\begin{equation}\label{colprod}
\text{if } c \in C(\sigma \mu), \psi(c)=\max_{a \in c} \varphi(a).
\end{equation}
\end{itemize}
We can now give the formulation of the main theorem of this article, which was conjectured by Stanley in \cite{St2} :
\begin{Th}\label{mainth}
Let $k,m$ be positive integers, $\mu \in S(k)$. For any sequences $\textbf{p},\textbf{q}$ of $m$ positive integers ($\textbf{q}$ being non-increasing), we have
\begin{equation}\label{thm}
\widehat{\chi}_{\textbf{p} \times \textbf{q}}(\mu) = (-1)^k \sum_{(\sigma,\varphi) \in S(k)^{(m)}} \left( \prod_{b \in C(\sigma)} p_{\varphi(b)} \prod_{c \in C(\sigma \mu)} -q_{\psi(c)} \right),
\end{equation}
where $\psi$ is defined by $(\ref{colprod})$.
\end{Th}
The case $m=1$ and the equality of highest degree terms (as polynomials in $\textbf{p}$ and $\textbf{q}$) have already been shown, respectively by R. Stanley in \cite{St1} and A. Rattan in \cite{Ra2}.\\
In section $\ref{sectii}$, we introduce some useful objects and results, which we will use in section $\ref{sectiii}$ to prove this theorem.

\section{Useful objects}\label{sectii}
\subsection{Stabilizers}
\indent The composition $\textbf{i} \circ \tau$ defines a right action of the
symmetric group on the sequences of integers between $1$ and $m$ of
length $k$. If $\textbf{i}$ is such a sequence, we will denote by $\Stab(\textbf{i})$
the stabilizer of $\textbf{i}$ (that is to say the set $\{\tau, \textbf{i} \circ
\tau=\textbf{i}\}$) and by $\delta_\textbf{i}$ its characteristic function :
$$\delta_\textbf{i}(\sigma)=\sum_{s \in \Stab(\textbf{i})} \delta_{\sigma s^{-1}},$$
where $\delta_\sigma=1$ if $\sigma=Id$ and $0$ else.
We will use the following easy properties :
\begin{itemize}
\item $\Stab(\textbf{i} \circ \tau)=\tau^{-1} \Stab(\textbf{i}) \tau$.
\item $|\{\tau,\textbf{i} \circ \tau=\textbf{j}\}|$ is either $0$ or $|\Stab(\textbf{i})|$.
\item Any sequence $\textbf{j}$ is in the orbit of exactly one non-increasing
  sequence.\\
\end{itemize}

\subsection{Young basis}
As usual, we draw a partition $\lambda \vdash k$ as a Young diagram ($\lambda_1$ squares in the first line, $\lambda_2$ in the second, and so on, all the lines are left justified). A Young standard tableau of shape $\lambda$ (their set will be denoted $YST(\lambda)$) is this Young diagram, filled with the numbers from $1$ to $k$, such that all lines and all columns are increasing. It is well-known that the dimension of the representation asso\-ciated to $\lambda$ is the cardinal of $YST(\lambda)$. We recall here the construction of a basis indexed by these objects :\\

Let $\lambda$ be a partition of $k$ and $W_\lambda$ be the irreducible associated $S(k)$-module (defined up to isomorphism). There exists on $W_\lambda$ an unique $S(k)$-invariant scalar product (the uniqueness is only true up to multiplication by a positive real number, but we will fix it for the end of the article). We will define by induction the Young orthonormal basis of $W_\lambda$. For $W_{(\textbf{0})}$, we have only one choice up to multiplication by a scalar. Then, we use the branching rule (see \cite{Sa}, theorem 2.8.3 for example) : as $S(k-1)$-module, we have
$$W_\lambda \simeq \bigoplus^\bot_{\substack{\lambda' \vdash k-1\\ \lambda' \leq \lambda}} W_\lambda',$$
where the inequality $\lambda' \leq \lambda$ is meant component by component (we will denote the two conditions by $\lambda' \nearrow \lambda$). The union of the Young orthonormal basis
for each $W_{\lambda'}$ is an orthogonal basis of $W_\lambda$. Multiplying each vector by a scalar, we obtain an orthonormal basis (we have a choice of unitary scalar to do but it doesn't matter). It is clear that the elements of the Young basis for $W_\lambda$ are indexed by sequences
$$\lambda_0=(0) \nearrow \lambda_1 \nearrow \ldots \nearrow \lambda_{k-1} \nearrow \lambda_k = \lambda,$$
or, equivalently, by Young standard tableaux of shape $\lambda$. We can denote this basis $(v_T)_{T \in YST(\lambda)}.$\\

As it is an orthonormal basis, we can use it to compute the character :
\begin{equation}\label{cartab}
\forall s \in \Q[S(k)], \chi_\lambda(s)=\tr_{W_\lambda}(s)=\sum_{T \in YST(\lambda)} \langle s \cdot v_T,v_T \rangle.
\end{equation}

\subsection{Jucys-Murphy elements}
The following elements of the symmetric group algebra, introduced by A. Jucys and G. Murphy (see \cite{Ju} and \cite{Mu}), play a very important role in the proof of the theorem :
$$\forall i \leq k,\text{ let }J_i = (1 i) + (2 i) + \ldots + (i-1 i).$$
These elements are very interesting because :
\begin{enumerate}
\item The products of different Jucys-Murphy elements have a nice combinatoric expression in the symmetric group algebra.
\item They are diagonal operators in the Young basis.
\end{enumerate}

\subsubsection{Combinatorics of products of Jucys-Murphy elements}
We will need in section $\ref{sectiii}$ the following result :

\begin{Lem}\label{jack}
For every integer $k$ and every set of variables $(X_j)_{1 \leq j\leq k}$, we have :
\begin{equation}\label{jackeq} X_1 \left(X_2 - J_2 \right) \ldots \left(X_k - J_k \right) = (-1)^k \sum_{\sigma \in S(k)} \left(\prod_{c \in C(\sigma)} -X_{min(c)} \right) \sigma.\end{equation}
\end{Lem}

Proof : we will prove it by induction over $k$, using the natural
embedding of $S(k-1)$ in $S(k)$. The case $k=1$ is obvious.\\
Let $k > 1$ and $\sigma \in S(k)$. We will look at the coefficient
of $\sigma$ in the left side of the equality $(\ref{jackeq})$. Using the induction
hypothesis, one can write :
\begin{eqnarray}
X_1 \left(X_2 - J_2 \right) \ldots \left(X_k - J_k
\right) \hspace{-3cm} & & \nonumber \\
&= & \left((-1)^{k-1} \sum_{\sigma' \in S(k-1)} \left(\prod_{c \in C(\sigma')}
-X_{min(c)} \right) \sigma' \right) \left(X_k - J_k
\right); \nonumber \\
& = & (-1)^{k} \sum_{\sigma' \in S(k-1)} \left(\prod_{c \in C(\sigma')}
-X_{min(c)} \right) (-X_k) \sigma' \nonumber \\
\label{combpr} && \qquad + (-1)^{k} \sum_{\sigma' \in S(k-1)}
\sum_{j=1}^{k-1} \left(\prod_{c \in C(\sigma')}-X_{min(c)}
\right) \sigma' (j k).
\end{eqnarray}
We distinguish two cases :
\begin{itemize}
\item $\sigma$ fixes $k$. It is the image of a permutation $\sigma' \in
  S(k-1)$ by the natural embedding. So, the coefficient of $\sigma$ in
  $(\ref{combpr})$ is
$$(-1)^{k} \left(\prod_{c \in C(\sigma')}
-X_{\min(c)} \right) \cdot (-X_k).$$
Looking at the cycles of the two permutations, we see that:
$$C(\sigma) = C(\sigma') \cup \{k\}.$$
So we have our result in this case.
\item Else, let $j=\sigma^{-1}(k)$, so that $\sigma$ can be written as $\sigma'
  (j k)$ and this is the only contribution to the coefficient of $\sigma$ in
  $(\ref{combpr})$. Thus $\sigma$ appears with the scalar :
$$(-1)^{k} \left(\prod_{c \in C(\sigma')} - X_{\min(c)} \right).$$
But the cycles of $\sigma'$ are the same as the cycles of $\sigma$, expected
that we have removed $k$ in the cycle which contained it (this cycle couldn't be the singleton $\{k\}$). So the
lemma is proved.

\end{itemize}

\subsubsection{Action on Young basis}
\begin{Def}
The content of the $j$-th box of the $i$-th line of a Young diagram is by definition the difference $j-i$. If $T$ is a standard tableau with $k$ boxes and $1 \leq a \leq k$, we will denote by $c_T(a)$ the content of the box containing the entry $a$.\\
\end{Def}
We can now state the following result, also due to A. Jucys and G. Murphy.
\begin{Lem}\label{murphy}
\begin{equation}\forall i \in \{1,\ldots,k\} \text{ and } T \in YST(\lambda), J_i(v_T) = c_T(i) v_T.\end{equation}
\end{Lem}
A proof can be found in [Ok2].\\

\emph{Remark :} an example of application of these two properties of Jucys-Murphy elements is the well-known formula :
$$\dim(\lambda) \prod_{\square \in \lambda} \left(X+c(\square)\right) = \sum_{s \in S(k)} \chi_\lambda(\sigma) X^{|C(\sigma)|}.$$

\subsection{Shifted Schur functions}
The other important objects in this paper are shifted Schur functions. The definition $\ref{defssf}$ and the theorem $\ref{formcar}$ can be found in A. Okounkov's and G. Olshanski's article \cite{OO} :
\begin{Defprop}\label{defssf}
Let $\lambda$ be a partition of $k$. We define the $\lambda$ shifted Schur polynomials as
$$s^\star_\lambda(x_1,x_2,\ldots,x_l)=\frac{\det \left( (x_i+l-i)_{\lambda_j +l-j} \right)_{1\leq i,j \leq l}}{\det \left( (x_i +l-i)_{l-j} \right)_{1\leq i,j \leq l}},$$
where $(r)_t = r(r-1)\ldots(r-t+1)$.
We can state that :
$$s^\star_\lambda(x_1,x_2,\ldots,x_l,0)=s^\star_\lambda(x_1,x_2,\ldots,x_l),$$
so $s^\star_\lambda$ is an element of the algebra of polynomials in countably many variables (symmetric in $x_1-1, x_2-2, \ldots$).
\end{Defprop}

\subsubsection{Link with characters}
As suggested by A. Rattan in his paper \cite{Ra1}, where he gives a new easier proof of the case $m=1$, we will use the following theorem :
\begin{Th}[A. Okounkov, G. Olshanski]\label{formcar}
For any integers $k \leq n$, permutation $\mu \in S(k)$ and partition $\nu \vdash n$, we have :
$$\widehat{\chi}_\nu(\mu)=\sum_{\lambda \vdash k} \chi_\lambda (\mu) s^\star_\lambda(\nu).$$
\end{Th}

\subsubsection{A new formula}
The new idea in this paper is to write $s^\star_\lambda(\nu)$ as the character of an element of the symmetric group algebra. We will be able to do this, thanks to this expression of shifted Schur polynomials, which was pointed to me by P. Biane and can be found in A. Okounkov's paper \cite{Ok1} (see (3.28) together with (3.34)) : for $\lambda \vdash k$,
\begin{eqnarray}
s^\star_\lambda(\nu)&=&\sum_{T \in YST(\lambda)} \sum_{i_1 \geq i_2 \geq \ldots \geq i_k \geq 1} \frac{1}{|\Stab(\textbf{i})|} \nonumber \\
\label{combssf} & & \qquad  \cdot \left[ \sum_{s \in \Stab(\textbf{i})} \langle s \cdot v_T, v_T  \rangle \nu_{i_1} \left(\nu_{i_2} - c_T(2) \right) \ldots \left(\nu_{i_k} - c_T(k) \right) \right].
\end{eqnarray}
Let, for $k \geq 1$ and $\nu \vdash n$ with $n \geq k$, $\mathcal{S}^k_\nu$ be the following element of ${\Q[S(n)]}$ :
\begin{equation}\label{relssf} \mathcal{S}^k_\nu = \sum_{i_1 \geq i_2 \geq \ldots \geq i_k \geq 1} \frac{\displaystyle \sum_{s \in \Stab(\textbf{i})} s}{|\Stab(\textbf{i})|} \nu_{i_1} \left(\nu_{i_2} - J_2 \right) \ldots \left(\nu_{i_k} - J_k \right),\end{equation}
Now, we can establish a new formula for shifted Schur polynomials.

\begin{Th}\label{newform}
If $\lambda \vdash k$ and $\nu \vdash n$,
\begin{equation} 
s^\star_\lambda(\nu)=\chi_\lambda(\mathcal{S}^k_\nu).
\end{equation}
\end{Th}
\emph{Proof :} Let's write the formula $(\ref{cartab})$ for $\mathcal{S}^k_\nu$.
$$\chi_\lambda(\mathcal{S}^k_\nu)=\sum_{T \in YST(\lambda)} \langle \mathcal{S}^k_\nu \cdot T, T \rangle.$$
The lemma $\ref{murphy}$ implies that, for any standard tableau $T$, any sequence $i_1 \geq i_2 \geq \ldots \geq i_k \geq 1$ and any $s \in \Stab(\textbf{i})$, we have
\begin{eqnarray*}
\left[ \nu_{i_1} \left(\nu_{i_2} - J_2 \right) \ldots \left(\nu_{i_k} - J_k \right) \right] \cdot v_T & \!\! = \!\! & \left[ \nu_{i_1} \left(\nu_{i_2} - c_T(2) \right) \ldots \left(\nu_{i_k} - c_T(k) \right) \right] v_T;\\
\left[s \cdot \nu_{i_1} \left(\nu_{i_2} - J_2 \right) \ldots \left(\nu_{i_k} - J_k \right) \right] \cdot v_T& \!\! = \!\! & \left[ \nu_{i_1} \left(\nu_{i_2} - c_T(2) \right) \ldots \left(\nu_{i_k} - c_T(k) \right) \right] s \cdot v_T.
\end{eqnarray*}
Now, taking the scalar product with $T$ and summing over $T$, $\textbf{i}$ and $s$, we find that $\chi_\lambda(\mathcal{S}^k_\nu)$ is exactly the right member of the equality $(\ref{combssf})$, so the theorem is proved.\\

The element $(\ref{relssf})$ might be very interesting to study to obtain results on shifted Schur functions. Here, we will look at the sum of the coefficients of permutations in the same conjugacy class as $\mu$.

\subsection{Orthogonality relations of the second kind}
In this paragraph, we will get from the orthogonality formula for irreducible characters an other relation. This can also be found in \cite{Sa}, thm 1.10.3. The classical formula can be written the following way :
for any partitions $\lambda$ and $\lambda'$ of $k$, we have
$$\sum_{\mathcal{C}} \frac{|\mathcal{C}|}{|G|} \chi_\lambda(\mathcal{C}) \chi_{\lambda'}(\mathcal{C})=\delta_{\lambda,\lambda'},$$
where the sum is taken over all conjugacy classes of $S(k)$. So we can reformulate it, saying that the square matrix $$\left(\chi_\lambda(\mathcal{C}) \sqrt{\frac{|\mathcal{C}|}{|G|}}\right)_{\substack{\lambda \vdash k \\ \mathcal{C} \text{ c.c. of }S(k)}}$$
is unitary. Looking at its rows, we have the following formula for any conjugacy classes $\mathcal{C}$ and $\mathcal{C}'$:
$$\sum_{\lambda \vdash k}  \frac{\sqrt{|\mathcal{C}| \cdot |\mathcal{C'}|}}{|G|} \chi_\lambda(\mathcal{C}) \chi_{\lambda}(\mathcal{C}')=\delta_{\mathcal{C},\mathcal{C}'}.$$
If $\mu$ and $\sigma$ are in the same conjugacy class $\mathcal{C}$, there are exactly $\frac{|G|}{|\mathcal{C}|}$ permutations $\tau \in S(k)$ such that $\tau \mu \tau^{-1} = \sigma$ (and of course there aren't any if they are in different conjugacy classes), so we can rewrite the previous equation under the form :

\begin{equation}
\sum_{\lambda \vdash k} \chi_\lambda(\mu) \chi_\lambda(\sigma)=\sum_{\tau \in S(k)} \delta_{\tau \mu \tau^{-1} \sigma^{-1}}.
\end{equation}

This formula can of course be extended by linearity in $\sigma$ to the group algebra $\Q[S(k)]$.

\begin{equation}\label{orth}
\sum_{\lambda \vdash k}  \chi_\lambda(\mu) \chi_\lambda \left( \sum_{\sigma \in S(k)} c_\sigma \sigma \right)=\sum_{\tau \in S(k)} \sum_{\sigma \in S(k)} c_\sigma \delta_{\tau \mu \tau^{-1} \sigma^{-1}}.
\end{equation}

\section{Proof of Stanley's conjecture}\label{sectiii}
As it is noticed in Stanley's paper \cite{St2}, one can reduce the conjecture to the case $p_1= \ldots =p_m=1$. Indeed, one can check easily that both sides of the equation $(\ref{thm})$ verify the following functional equation in $\textbf{p}$ and $\textbf{q}$ :
$$F(\textbf{p},\textbf{q})|_{q_i=q_{i+1}} = F(p_1,\ldots,p_{i}+p_{i+1},\ldots,p_m,q_1,\ldots,q_{i},q_{i+2},\ldots,q_m).$$
In this case, the partition $\nu=1^m \times \textbf{q}$ is simply given by $\nu_i=q_i$ if $1 \leq i \leq m$ and $\nu_i=0$ else.\\

As mentionned earlier, we will use theorems $\ref{formcar}$ and $\ref{newform}$ to get the following expression of the normalised character :
\begin{equation}\widehat{\chi}_{1^m \times \textbf{q}}(\mu) = \sum_{\lambda \vdash k}
\chi_\lambda(\mu) \chi_\lambda (\mathcal{S}^k_{1^m \times \textbf{q}}).\end{equation}
The lemma $\ref{jack}$, applied to $X_j=q_{i_j}$ gives a nicer expression of $\mathcal{S}^k_{1^m \times \textbf{q}}$ :
\begin{eqnarray}\widehat{\chi}_{1^m \times \textbf{q}}(\mu) & = &  \\
& & \hspace{-2,2 cm}(-1)^k \sum_{\lambda \vdash k} \chi_\lambda(\mu)
  \chi_\lambda\left(\sum_{m \geq i_1 \geq i_2 \geq \ldots \geq i_k \geq 1}
  \frac{\displaystyle \sum_{s \in \Stab(\textbf{i})} s}{|\Stab(\textbf{i})|}
  \sum_{\sigma \in S(k)} \left(\prod_{c \in C(\sigma)} -q_{i_{\min(c)}}
  \right) \sigma \right). \nonumber \end{eqnarray}
We now use the orthogonality relation $(\ref{orth})$ and obtain :
\begin{eqnarray} \widehat{\chi}_{1^m \times \textbf{q}}(\mu) & = & (-1)^k \sum_{\tau \in S(k)}\sum_{m \geq i_1 \geq i_2 \geq \ldots
  \geq i_k \geq 1} \frac{1}{|\Stab(\textbf{i})|} \cdot \nonumber \\ & & \qquad \left[ \sum_{\sigma \in S(k)} \left(\prod_{c \in C(\sigma)}
  -q_{\max\limits_c \textbf{i}}\right) \sum_{s \in \Stab(\textbf{i})} \delta_{\tau \mu \tau^{-1} \sigma^{-1} s^{-1}} \right]
\end{eqnarray}
As $\sigma \mapsto \tau \sigma \tau^{-1}$ defines a bijection of the symmetric group into itself, we can replace $\sigma$ by $\tau \sigma \tau^{-1}$ in the second line of the previous equation.
\begin{eqnarray}
\widehat{\chi}_{1^m \times \textbf{q}}(\mu) & = & (-1)^k \sum_{\tau \in S(k)}\sum_{m \geq i_1 \geq i_2 \geq \ldots
  \geq i_k \geq 1} \frac{1}{|\Stab(\textbf{i})|} \cdot \nonumber \\
&  & \qquad \left[ \sum_{\sigma \in S(k)}
  \left(\prod_{c \in C(\tau \sigma \tau^{-1})}
  -q_{\max\limits_c \textbf{i}}\right) \delta_{\textbf{i}}(\tau \mu \sigma^{-1} \tau^{-1}) \right]\\
\widehat{\chi}_{1^m \times \textbf{q}}(\mu) & = & (-1)^k \sum_{\tau \in S(k)}\sum_{m \geq i_1 \geq i_2 \geq \ldots
  \geq i_k \geq 1} \nonumber \\& & \quad \left[\frac{1}{|\Stab(\textbf{i} \circ \tau)|} \sum_{\sigma \in S(k)}
  \left(\prod_{c \in C(\sigma)}
  -q_{\max\limits_c \textbf{i} \circ \tau}\right) \delta_{\textbf{i} \circ \tau}(\mu \sigma^{-1})\right]\end{eqnarray}
Now, we can see that the expression in the brackets depends only on the
  sequence $\textbf{i} \circ \tau$. Each sequence $(j_1,\ldots,j_k)$ of integers
  between $1$ and $m$ (but not necessarily non-increasing) can be
  written as $|\Stab(\textbf{j})|$ different ways as $\textbf{i} \circ \tau$, where $\textbf{i}$
  is non-increasing and $\tau \in S(k)$ (in all these writings $\textbf{i}$ is the same but $\tau$ can be chosen among $|\Stab(\textbf{j})|=
|\Stab(\textbf{i})|$ permutations). So we have :
\begin{eqnarray}\widehat{\chi}_{1^m \times \textbf{q}}(\mu) & = & (-1)^k \sum_{\substack{1
      \leq j_1 \leq m\\ \vdots\\ 1 \leq j_k \leq m}} \sum_{\sigma \in S(k)}
      \left(\prod_{c \in C(\sigma)}  -q_{\max\limits_c \textbf{j}}\right)
      \delta_{\textbf{j}}(\mu \sigma^{-1});\\
\widehat{\chi}_{1^m \times \textbf{q}}(\mu) & = & (-1)^k \sum_{\sigma \in S(k)} \sum_{\textbf{j} \text{ fixed by } \mu
	\sigma^{-1}} \left(\prod_{c \in C(\sigma)}  -q_{\max\limits_c \textbf{j}}\right).\end{eqnarray}
Note that the sequences of integers between $1$
      and $m$ fixed by a permutation are
      exactly the colourings of its cycles in $m$ colours (the
      colour of a cycle is the common value of $\textbf{j}$ on its elements). So, if we change the index of the sum over $S(k)$ putting $\sigma'=\sigma \mu^{-1}$, we can write :

\begin{eqnarray}\widehat{\chi}_{1^m \times \textbf{q}}(\mu) & = & (-1)^k \sum_{\sigma' \in S(k)} \sum_{\substack{\textbf{j}\text{ such that}\\(\sigma',\textbf{j}) \ \in \ S(k)^{(m)}}} 
      \left(\prod_{c \in C(\sigma' \mu)}  -q_{\max\limits_c \textbf{j}}\right),\end{eqnarray}
which is exactly $(\ref{thm})$ in the case $p_1=\ldots=p_m=1$.

\end{document}